\documentstyle[12pt]{article}

\begin{document}
\newcommand{\ol }{\overline}
\newcommand{\ul }{\underline }
\newcommand{\ra }{\rightarrow }
\newcommand{\lra }{\longrightarrow }
\newcommand{\ga }{\gamma }
\newcommand{\st }{\stackrel }
\title{{\Large\bf Some Notes on the Baer-invariant of a Nilpotent Product of
Groups}\footnote{This research was in part supported by a grant
from IPM}}
\author{Behrooz Mashayekhy\\ Department of Mathematics, Ferdowsi University of
Mashhad,\\ P.O.Box 1159-91775, Mashhad, Iran\\E-mail: mashaf@science2.um.ac.ir\\
and \\ Institute for Studies in Theoretical Physics and Mathematics,\\
Tehran, Iran }
\date{ }

\maketitle

\begin{abstract}
 W.Haebich (Bull. Austral. Math. Soc., 7, 1972, 279-296) presented a formula for
the Schur multiplier of a regular product of groups. In this
paper first, it is shown that the Baer-invariant of a nilpotent
product of groups with respect to the variety of nilpotent groups
has a homomorphic image and in finite case a subgroup of
Haebich's type. Second a formula will be presented for the
Baer-invariant of a nilpotent product of cyclic groups with
respect to the variety of nilpotent groups.\\
A.M.S. Classification: 20F12,20F18,20K25 
\end{abstract}
\begin{center}
{\bf 1.INTRODUCTION}
\end{center}

 In 1907, I.Schur [ ], using representation method, showed that if $G$ is the
direct product of $A$ and $B$, then the Schur multiplier of
$G$ has the following form :  \\
   $$ M(G)=M(A\times B)\cong M(A)\times M(B)\times (A_{ab}\otimes B_{ab}).$$
Also, J.Wiegold [ ] in 1971 obtained the above result by some properties of
covering groups:
$$ M(A\times B)\cong M(A)\oplus M(B)\oplus \frac {[A,B]}{[A,B,A*B]}\ ,\
where\ \frac {[A,B]}{[A,B,A*B]}\cong A_{ab}\otimes B_{ab}\ . $$

 In 1979, M.R.R.Moghaddam [ ] found a formula for the Baer-invariant ( the
generalization of the Schur multiplier with respect to the variety of groups.)
of a direct product of two groups with respect to the variety of nilpotent
groups of class at most $c$, ${\cal N}_c$ , where $c+1$ is a prime number or
$4$. This result generalized somehow the work of Schur and Wiegold.

 In 1997, we [ ] presented an explicit formula for the Baer-invariant of a
finite abelian group with respect to the variety of nilpotent groups for every
$c\geq 1$ . Since a finite abelian group is a direct product of some cyclic
groups, so this result is somehow a generalization of Schur-wiegold's result and
also Moghaddam's result.

 It is known that nilpotent product is a generalization of direct product.
Therefore it is interesting to find a formula for the Schur multiplier or
the Baer-invariant of a nilpotent product.

 In 1972, W.Haebich [ ] (Theorem 2.14) presented a formula for the Schur
multiplier of a regular product of a family of groups. We know that regular
product is a generalization of nilpotent product, so this result is an
interesting generalization of the Schur-Wiegold's result. Also, M.R.R.Moghaddam
[ ] (Theorem 2.7) in 1979 gave a formula similar to Haebich's formula for
the Schur multiplier of a nilpotent product. His approach was different from
W.Haebich.

 Finally, N.D.Gupta and M.R.R.Moghaddam [ ] (Theorem 2.15) in 1992 have
calculated the Baer-invariant of the nilpotent dihedral group of class $n$
$$ G_n=<x,y,|x^2,y^2,[x,y]^{2^{n-1}}> $$
with respect to the variety of nilpotent groups, ${\cal N}_c$. It is routine to
verify that $G_n$ is isomorphic to the $n$th nilpotent product of two cyclic
groups i.e $G_n={\bf Z}_2\st {n}{*} {\bf Z}_2$.\\
\ \ \\
{\bf 2.NOTATION AND PRELIMINERIES}

 It is supposed that the reader is familiar with the notions of variety of
groups, verbal and marginal subgroups, Schur multiplier and Baer-invariant. In
particular if $G$ is a group with a free presentation
$$ 1\lra R\lra F\lra G\lra 1\ , $$
then the Baer-invariant of $G$ with respect to the variety of nilpotent groups
of class at most $c\geq 1\ ,{\cal N}_c$, denoted by ${cal N}_cM(G)$, is defined
to be
$$ {\cal N}_cM(G)=\frac {R\cap \ga_{c+1}(F)}{[R,\ _cF]}  $$
(See C.R.Leedham-Green and S.McKay [ ], from where our notation have been taken
and H.Neumann [ ] for the notion of variety of groups.)

 Clearly, if $c=1$, then ${\cal N}_c=\cal A$ is the variety of abelian groups
and
$$ {\cal N}_cM(G)=M(G)\cong \frac {R\cap F'}{[R,F]} $$
is the {\it Schur multiplier} of $G$.\\
{\bf Definition 2.1}

 A variety $\cal V$ is said to be a {\it Schur-Baer variety} if for any group
$G$ for which the marginal factor group, $G/V^*(G)$, is finite, then the verbal
subgroup,$V(G)$, is also finite and $|V(G)|$ divides a power of $|G/V^*(G)|$.
I.Schur in [11] proved that the variety of abelian groups,$\cal A$, is a
Schur-Baer variety. Also R.Baer in [1] proved that the variety defined by some
outer commutator words has the above property.
 The following theorem tells us a very important property of Schur-Baer
varieties.\\
{\bf Theorem 2.2}(C.R.Leedham-Green and S.McKay [6])

 The following conditions on the variety $\cal V$ are equivalent:\\
$(i)$ $\cal V$ is a Schur-Baer variety.\\
$(ii)$ For every finite group $G$, its Baer-invariant, ${\cal V}M(G)$, is of
order dividing a power of $|G|$\ .\\
{\bf Definition 2.3}

Consider the following map
$$ \varphi^*:\prod_{i\in I}^{*} A_i\lra \prod_{i\in I}^{}\!^{\times }A_i$$
$$ a_1a_2\ldots a_n\longmapsto (a_1,a_2,\ldots ,a_n) $$
which is a natural map from the free product of $\{A_i\}_{i\in I}$ onto the
direct product of $\{A_i\}_{i\in I}$. Clearly its kernel is the normal closure
of
     $$ <[A_i,A_j]\ |\ i,j\in I,i\neq j> $$
in the free product $\prod_{i\in I}^{*}\!A_i$ and denoted by $[A_i]^*$ and
called {\it the cartesian subgroup} of the free product.

 For the properties of cartesian subgroup see H.Neumann [ ].\\
{\bf Definition 2.4}

 Let $\cal V$ be a variety of groups defined by a set of laws $V$ then {\it the
verbal product} of a family of groups $\{A_i\}_{i\in I}$ associated with the
variety $\cal V$ is defined to be
$$ {\cal V}\prod_{i\in I}^{}A_i = \frac {\prod_{i}^{*}A_i}{V(A)\cap [A_i]^*}\ . $$

 The verbal product is also known as {\it varietal product} or simply $\cal
V$\_{\it product}. If $\cal V$ is the variety of all groups, then the
correspoding verbal product is the {\it free product}; if $\cal V=A$ is the
variety of all abelian groups, then the verbal product is the {\it direct
product}. We denote ${\cal V}\prod_{i}^{}A_i$ for the  $\cal V$\_product of the
family of groups $A_i$, $i\in I$.

 The $\cal V$-product of the ${A_i}'s$ is `between' the free product and the
direct product in the sense that the natural epimorphism of the
free product onto the direct product can be factored through the
$\cal V$\_product, i.e
$$ \prod_{i}^{ }\!^*A_i \st {\nu }{\lra } {\cal V}\prod_{i}^{}A_i \st {\varphi
(\nu )}{\lra } \prod_{i}^{}\!^{\times }A_i\ \ . $$
The kernel of the natural epimorphism $\varphi (\nu )$ of the $\cal V$\_product
onto the direct product of the ${A_i}'s$ is the cartesian subgroup of ${\cal
V}\prod_{i}^{}A_i$, that is the normal closure of
$<[A_i,A_j]\ |\ i\neq j,i,j\in I>$, in the $\cal V$-product, which is denoted by
$[A_i]^V$\ .\\
{\bf Lemma 2.5}

 With the above notation, we have\\
$(i)$ $Ker\varphi(\nu)=[A_i]^V=\nu ([A_i]^*)$ where $\varphi^*=\varphi (\nu
)\circ \nu $ as given in 2.7 .\\
$(ii)$ If $a\neq 1$ is an element of ${\cal V}\prod_{i}^{}A_i$, then
$$ a=a_{i_1}\ldots a_{i_m}c\ \ ,\ \ 1\neq a_{i_j}\in A_{i_j}\ \ ,i_1<i_2<\ldots
<i_m\ \ ,\ \ c\in [A_i]^V\ , $$
and the $a_{i_j}$ and $c$ are uniquely determined by $a$ and chosen order of
$I$\ .\\
{\bf Proof.} See [ ] .\\
{\bf Definition 2.6}

 Let $\{A_i\}_{i\in I}$ be a family of groups and ${\cal N}_c$ be the variety
of nilpotent groups of class at most $c\geq 1$. Then the ${\cal N}_c$-product,
${\cal N}_c\prod_{i}^{}A_i$, is called the $c${\it th nilpotent product} of
$A_i$'s. In particular, if $A$ and $B$ are two groups, then the $c$th
nilpotent product of $A$ and $B$, denoted by $A\st {c}{*}B$, is as follows:
$$ A\st {c}{*}B=\frac {A*B}{[A,B,\ _{c-1}A*B]}\ . $$
{\bf Theorem 2.7} (M.R.R.Moghaddam [9])

 Let $A$ and $B$ be two groups and $A\st {n}{*}B$ be the $n$th nilpotent
product of $A$ and $B$. Then for $n\geq 1$
$$ M(A\st {n}{*}B)\cong M(A)\oplus M(B)\oplus \frac {[A,B,\
_{(n-1)}A*B]}{[A,B,\ _nA*B]}\ . $$
{\bf Definition 2.8}

 Let $\{A_i\}_{i\in I}$ be a family of subgroups of an arbitrary group $G$. We
say that $G$ is a {\it regular product} of its subgroups ${A_i}'s,i\in I$, where
$I$ is an ordered set if the following two conditions hold:\\
$(i)\ \ G=<A_i|i\in I>,$\\
$(ii)\ \ A_i\cap \hat {A_i}=1$\ \ for all\ \ $i\in I$,\ \ where\ \ $\hat
{A_i}=\prod_{j\neq i,j\in I}A_j^G$\\
The subgroups ${A_i}'s,i\in I$ will be referred to as the {\it regular factors} of
the group $G$. The direct product, the free product, and the verbal product of
an arbitrary set of groups are examples of regular products.\\
{\bf Lemma 2.9}

Let $A$ be a subgroup of a group $G$, and $N$ be a normal subgroup of $G$ and
$\{M_i|i\in I\}$ be a family of normal subgroups of $G$. Then
$$ [A\prod_{i}{}M_i,N]=[A,N]\prod_{i}{}[M_i,N]\ .$$

 The following theorem gives a characterization of regular product .\\
{\bf Theorem 2.10} (O.N.Golovin [2])

 Suppose that a group $G$ is generated by a family $\{A_i|i\in I\}$ of its
subgroups, where $I$ is an ordered set. Then $G$ is a regular product of the
${A_i}'s$ if and only if every element of $G$ can be written uniquely as a product
$$ a_{\lambda_1}a_{\lambda_2}\ldots a_{\lambda_m}u\ ,$$
where $$1\neq a_{\lambda_i}\in A_{\lambda_i},\lambda_i< \lambda_2< \ldots
<\lambda_m $$
and $$u\in [A_i^G]=<[A_{\lambda }^G,A_{\mu }^G]|\lambda \neq \mu ,\lambda
,\mu \in I>\ .$$\\
{\bf Proof.} See [2,5] .\\
{\bf Definition 2.11}

 Let $G$ be a regular product of ${A_i}'s,i\in I$. A homomorphism $f:G\lra H$ of
$G$ to a group $H$ is called a {\it regular homomorphism} if
  $$Ker\,f\subseteq [A_i^G]\ .$$
The groups $Im\,f$ and $G/Ker\,f$ are called a {\it regular homomorphic image}
and a {\it regular quotient group}, respectively.\\
{\bf Theorem 2.12}

 If $G$ is a regular product of the subgroups ${A_i}'s,i\in I$, and $f:G\lra H$ is
a regular homomorphism, then $f$ restricted to $A_i$ is an {\it isomorphism}
for each $i\in I$ and $f(G)$ is a regular product of $f(A_i)$.\\
{\bf Proof.} See [10] .\\
{\bf Theorem 2.13}

 Let $G$ be generated by the ${A_i}'s,i\in I$, and let
$$ \psi :\prod_{i\in I}^{ }\!^{*}A_i\lra G$$
be the homomorphism induced by the identity map on each $A_i$. Then $G$ is a
regular product of the $A_i$'s if and only if $\psi $ is regular.\\
{\bf Proof.} See [ ].\\
{\bf Theorem 2.14} (W.Haebich [ ])

 Let $G$ be a regular product of a family of its subgroups $\{A_i|i\in I\}$
and let $H$ be the kernel of the natural homomorphism
$$ A=\prod_{i}^{ }\!^{*}A_i\lra G\ \ \ ,$$
induced by the identity map on each $A_i\ ,i\in I$\ . Then
$$ M(G)\cong (\prod_{i\in I}^{ }\!^{\times }M(A_i))\oplus \frac {H}{[H,A]}\ .$$
{\bf Theorem 2.15}(N.D.Gupta and M.R.R.Moghaddam [ ])

 Let $G_n$ be the nilpotent dihedral group of class $n$, i.e $G_n={\bf Z}_2\st
{n}{*} {\bf Z}_2$, then
$(i)$ If $c<n$, then ${\cal N}_cM(G_n)$ is the cyclic group of order $2^c$.
$(ii)$ If $c\geq n$, then ${\cal N}_cM(G_n)$ is an abelian extension of a
cyclic group of order $2^{n-1}$ by an elementary abelian 2-group of rank
$r(c+1)$, the number of basic commutators of wieght $c+1$ on two letters.
In particular, for $c\geq n$ ,
$$ {\cal N}_cM(G_n)\cong {\bf Z}_{2^n}\oplus {\bf Z}_2\oplus \ldots \oplus
{\bf Z}_2\ (\ r(c+1)-1\ copies\ of\ {\bf Z}_2\ )\ .$$
{\bf Proof.} See [ ] Theorem   . \\
{\bf Lemma 2.16}

If $H$ is a subgroup of a finite abelian group $G$, then $G$ has a subgroup
isomorphic to $G/H$ .\\
\ \ \ \\
{\bf 3. A SUBGROUP OF HAEBICH'S TYPE}

 Let $\{A_i|i\in I\}$ be a family of groups (we can consider $I$ as an ordered
set) and for each $i\in I$, $F_i$ denote a fixed free group such that the
following exact sequence be a free presentation for $A_i$
$$ 1\lra R_i\lra F_i\st {\nu_i}{\lra} A_i\lra 1\ \ . $$
We denote by $\nu$ the natural homomorphism from the free product
$F=\prod_{i\in I}^*F_i$ onto $A=\prod_{i\in I}^*A_i$ induced by the
$\nu_i$'s.Also,
the group $G_n$ will be assumed to be the $n$th nilpotent product of the
$A_i$'s, $i\in I$. If $\psi_n$ is the natural homomorphism from $A$ onto $G_n$
induced by the identity map on each $A_i$ ;
$$ F=\prod_{i\in I}^*F_i\st {\nu }{\lra }A=\prod_{i\in I}^*A_i\st {\psi_n}{\lra}
G_n=\prod_{i\in I}^{\st {n}{*}}A_i\lra 1\ ,$$
then we denote by $H_n$ be the kernel of $\psi_n$ and $L_n$ the inverse image of
$H_n$ in $F$ under $\nu$ ; $ H_n=ker\psi_n \ \ \& \ \ L_n=\nu^{-1}(H_n)\ .$

Now we have the following lemma.\\
{\bf Lemma 3.1}

 By the above notation and assumption, we have\\
$(i)\ H_n=\nu (K_n)$ , where $K_n=\ga_{n+1}(F)\cap [F_i]^*$\\
$(ii)\ G_n\cong F/L_n$ , where $L_n=ker(\psi_n \circ \nu)=(\prod_{i\in I}
R_i^F)K_n $\\
{\bf Proof.}\\
$(i)$ By definition of $n$th nilpotent product, we have
$H_n=\ga_{n+1}(A)\cap [A_i]^*$, where $[A_i]^*$ is the cartesian subgroup of
$A$ , so clearly $\nu (K_n)=H_n$ .\\
$(ii)$ Clearly $Ker\nu =\prod_{i\in I} R_i^F$ , and composition
$\psi_n\circ \nu:F\lra G_n$ is surjective homomorphism. Since $H_n=Ker\psi_n
=\nu (K_n)$ , it follows that
$$ L_n=Ker(\psi_n\circ \nu )=(Ker\nu )K_n=(\prod_{i\in I} R_i^F)K_n\ \ .$$
Hence $G_n\cong F/L_n$ .\ \ $\Box$\\
{\bf Lemma 3.2}

 If $ D_1=\prod_{i \neq j\ ,\ i,j \in I}[R_i,F_j]^F$ , then $\prod_{i\in
I}R_i^F=(\prod_{i\in I}R_i)D_1$ .\\
{\bf Proof.}

 See [ ] Lemma   .\\
{\bf Notation 3.3}

 We define
$(i)\  D_c=\prod_{\exists j,\mu_j\neq i}[R_i,F_{\mu_1},\ldots ,F_{\mu_c}]^F\
,$\\
$(ii)\ E_c=D_1\cap \ga_{c+1}(F) \ .$\\
Clearly $D_c\subseteq E_c$ for every $c\geq 1$ .

 Now we prove the following important lemma.\\
{\bf Lemma 3.4}

 By the above notation and assumotion, we have\\
$(i)\ [L_n,\ _cF]=(\prod_{i\in I}[R_i,\ _cF_i])D_c[K_n,\ _cF]\ .$ \\
$(ii)$ If $n\leq c$ , then $L_n\cap \ga_{c+1}(F)=(\prod_{i\in I}(R_i\cap
\ga_{c+1}(F_i))E_cK_c$ .\\
$(iii)$ If $n>c$ , then $L_n\cap \ga_{c+1}(F)=(\prod_{i\in
I}(R_i\cap \ga_{c+1}(F_i))E_cK_n$ .\\
{\bf Proof.}

 $(i)$ By definition of $C_{(n)}$ , we have
$$ [C_{(n)},\ _cB]=[(\prod_{i\in I}C_i^B)K_n,\ _cB]=[(\prod_{i\in I}C_i^B),\
_cB][K_n,\ _cB]\ \ \ \ \ \ ({\it by\ Lemma\ 2.18})$$
$$ \ \ \ \ \ \ \ \ \ =(\prod_{i\in I}[C_i,\ _cB]^B)[K_n,\ _cB]\ \ \ \ \ \
({\it by\ Lemma\ 2.18})$$
$$ \ \ \ \ \ \ \ \ \ \ \ \ \ \subseteq (\prod_{i\in I}[C_i,\ _cB_i])D_c[K_n,\
_cB]\ \ \ \ \ ({\it by\ definition\ of\ D_c})$$
The reverse inclusion is obviously true, so the result holds.

$(ii)$ Let $n\leq c$ , and $g\in C_{(n)}\cap \ga_{c+1}(B)$ , by Lemmas 3.1 and
3.2 we have
$$ g=c_{\lambda_1}\ldots c_{\lambda_t}dk\ \ ,$$
where $c_{\lambda_i}\in C_{\lambda_i}\ ,\ d\in D_1\ ,\ k\in K_n$ . Now consider
the following homomorphism
$$ B=\prod_{i\in I}^{ }\!^*B_i \st{nat}{\lra} \prod_{i\in I}^{ }\!^{\times}B_i
\ \ .$$
Since $g\in \ga_{c+1}(B)$ , the image of $g$ under the above natural
homomorphism is
$$ (c_{\lambda_1},\ldots ,c_{\lambda_t})\in \ga_{c+1}(\prod
\!^{\times}B_i)=\prod \!^{\times}\ga_{c+1}(B_i)\ \ .$$
Therefore $c_{\lambda_i}\in \ga_{c+1}(B_{\lambda_i})\cap C_{\lambda_i}$ and so
$dk\in \ga_{c+1}(B)\cap [B_i]^*=K_c$ .\\
The reverse inclusion can be seen easily, and so the result holds.

$(iii)$ It is clear.\ \ $\Box$\\


\begin{thebibliography}{14}

\bibitem{1} R.Baer:``Endlichkeitskriterien f\"{u}r Kommutatorgruppen.'' Math.
Ann., 124,161-177 (1952).

\bibitem{2} O.N.Golovin:`` On the Isomorphisms of Nilpotent Decompositions of
Groups.'' Mat. Sbornik N.S. 28(70),445-52 (1951) Amer. Math. Soc.
Translations,Ser. 2,2,133-140 (1956).

\bibitem{3} W.Haebich:``The Multiplicator of a Regular Product of Groups.''
Bull. Austral. Math. Soc.,7,279-296,(1972).

\bibitem{4} T.W.Hungerford:``Algebra'' Springer-Verlag,G.T.M. 73 (1980).

\bibitem{5} G.Karpilovsky:``The Schur Multiplier.'' London Math. Soc.
Monographs,New Series no. 2,(1987).

\bibitem{6} C.R.Leedham-Green and S.McKay:``Baer-invariant,Isologism,
Varietal Laws and Homology.'' Acta Math.,137,99-150 (1976).

\bibitem{7} B.Mashayekhy and M.R.R.Moghaddam:``Higher Schur Multiplicator of
a Finite Abelian Group'',Algebra Colloquium, 4:3, 317-322, 1997.

\bibitem{8} M.R.R.Moghaddam:``The Baer-invariant of a Direct Product.'' Arch.
der. Math.,33,504-511,(1979).

\bibitem{9} M.R.R.Moghaddam:``The Schur Multiplicator of Nilpotent Products.''
Bull. Iranian math. Soc. ,vol. 6,no. 2,79-89,(1979).

\bibitem{10} H.Neumann:``Varieties of Groups.'' Berlin,Springer-Verlag,(1967).

\bibitem{11} I.Schur:``\"{U}ber die Darstellung der Endlichen Gruppen durch
Gebrochene Lineare Substitutionen.'' J. Reine Angew. Math.,127,20-50,(1904).

\bibitem{12} I.Schur:``Untersuchungen \"{u}ber die Darstellung der Endlichen
Gruppen durch Gebrochene Lineare Substitutionen.'' J. Reine Angew. Math.,132,
85-137,(1907).

\bibitem{13} R.F.Turner-Smith:``Marginal Subgroup Propeties for Outer
Commutator Words.'' Proc. London Math. Soc.,14,321-324,(1964).

\bibitem{36} J.Wiegold:``The Multiplicator of a Direct Product.'' Quart. J.
Math. Oxford,(2),22,103-105,(1971).
\end{thebibliography}
\end{document}